\documentclass[11pt]{article}
\usepackage{mathrsfs}
\usepackage{nicefrac}
\usepackage[T1]{fontenc}
\usepackage{latexsym,amssymb,amsmath,amsfonts,amsthm}
\usepackage{graphicx}
\usepackage{caption}
\usepackage[utf8]{inputenc}
\usepackage{authblk}
\begin{document}

\centerline{\bf\Large An improved upper bound for the number of distinct eigenvalues}
\centerline{\bf\Large of a matrix after perturbation}

\vskip 0.3cm

\centerline{Xuefeng Xu}
\centerline{\small Institute of Computational  Mathematics and Scientific/Engineering Computing,}
\centerline{\small Academy of Mathematics and Systems Science, Chinese Academy of Sciences,}
\centerline{\small  Beijing 100190, China. E-mail address: xuxuefeng@lsec.cc.ac.cn.}

\vskip 0.3cm

\noindent{\bf \large Abstract}

An upper bound for the number of distinct eigenvalues of a perturbed matrix has been recently established by P. E. Farrell [1, Theorem 1.3]. The estimate is the central result in Farrell’s work and can be applied to estimate the number of Krylov iterations required for solving a perturbed linear system. In this paper, we present an improved upper bound for the number of distinct eigenvalues of a matrix after perturbation. Furthermore, some results based on the improved estimate are presented.

\vskip 0.2cm

\noindent{\bf Keywords:} Distinct eigenvalues; Perturbation;  Defectivity; Derogatory index

\vskip 1.2cm

\noindent{\bf \large 1.\ Introduction}

\vskip 0.2cm

The spectrum of a matrix after perturbation has been investigated by many authors. However, most work is devoted to discussing some special cases, especially the case of symmetric rank-one perturbations; see, for instance, [2-5]. Recently, P. E. Farrell [1] presented an upper bound for the number of distinct eigenvalues of arbitrary matrices perturbed by updates of arbitrary rank. Let $\mathbb{C}^{n\times n}$, $\Lambda(\cdot)$, ${\rm rank}(\cdot)$, and $|\cdot|$ be the set of all $n\times n$ complex matrices, the set of all distinct eigenvalues of a square matrix, the rank of a matrix, and the cardinality of a set, respectively.
Assume that $A,B\in\mathbb{C}^{n\times n}$, and let $C=A+B$, it is proved by Farrell [1, Theorem 1.3] that
\begin{align*}
|\Lambda(C)|\leq({\rm rank}(B)+1)|\Lambda(A)|+d(A),\tag{1.1}
\end{align*}
where $d(\cdot)$ denotes the defectivity of a matrix (see Definition 2.2 below). The result can be used to estimate the number of Krylov iterations for solving a linear system. 

It follows from the definitions of $\Lambda(\cdot)$ and $|\cdot|$ that $|\Lambda(M)|\leq n$ for all $M\in\mathbb{C}^{n\times n}$. Given $A\in\mathbb{C}^{n\times n}$, we can observe that the estimate (1.1) is mainly of interest in the situation that ${\rm rank}(B)$ is small, that is to say, $B$ is a low-rank perturbation. More specifically, if ${\rm rank}(B)\leq\frac{n-d(A)}{|\Lambda(A)|}-1$, then $({\rm rank}(B)+1)|\Lambda(A)|+d(A)$ $(\leq n)$ is an applicable upper bound. On the other hand, if ${\rm rank}(B)>\frac{n-d(A)}{|\Lambda(A)|}-1$, then $({\rm rank}(B)+1)|\Lambda(A)|+d(A)$ $(> n)$ is a trivial upper bound.

Nevertheless, the estimate (1.1) is not sharp in certain cases. We now give a specific example to illustrate the defect of the upper bound in (1.1). We choose a matrix $A$ as follows:
\[A=\begin{pmatrix}
\lambda_{0} & 1 & 0 & \cdots & 0\\
0 & \lambda_{0} & 1 & \cdots & 0\\
\vdots & \vdots & \ddots & \ddots & \vdots \\
0 & 0 & \cdots & \lambda_{0} & 1\\
0 & 0 & \cdots & 0 & \lambda_{0}\\
\end{pmatrix},\]
which is a $n\times n$ Jordan block. Thus, $|\Lambda(A)|=1$ and $d(A)=n-1$. Let $n\times n$ matrix $B_{r}$ with rank $r$ $(1\leq r\leq n-1)$ be defined by
\[B_{r}=\begin{pmatrix}
T_{r} & O_{r\times (n-r-1)}\\
O_{(n-r)\times(r+1)} & O_{(n-r)\times (n-r-1)}\\
\end{pmatrix},\ \ \ \ T_{r}=\begin{pmatrix}
1 & -1 & 0 & \cdots & 0\\
0 & 2 & -1  & \cdots & 0\\
\vdots & \vdots & \ddots & \ddots & \vdots \\
0  & 0 & \cdots & r  & -1 \\ 
\end{pmatrix}.\]
Hence, the upper bound in the estimate (1.1) is $({\rm rank}(B_{r})+1)|\Lambda(A)|+d(A)=n+r>n$. In this case, the upper bound in (1.1) is always invalid (i.e., the upper bound is strictly greater than order $n$) for all $1\leq r\leq n-1$.

In this paper, we give an improved upper bound for the number of distinct eigenvalues of a matrix after perturbation. Under the same assumptions, we establish that
\begin{align*}
|\Lambda(C)|\leq({\rm rank}(B)+1)|\Lambda(A)|+d(A)-d(C).\tag{1.2}
\end{align*}
Applying (1.2) to the above example, we can derive that the improved upper bound of $|\Lambda(C_{r})|$ (here $C_{r}=A+B_{r}$) is $({\rm rank}(B_{r})+1)|\Lambda(A)|+d(A)-d(C_{r})=(n+r)-(n-1-r)=2r+1$, which is an applicable upper bound, especially in low-rank perturbations. 

We give another numerical example to explain our estimate in (1.2). Let
\[A=\begin{pmatrix}
1 & 1 & 0 & 0 & 0\\
0 & 1 & 1 & 0 & 0\\
0 & 0 & 1 & 0 & 0\\
0 & 0 & 0 & 1 & 0\\
0 & 0 & 0 & 0 & 1\\
\end{pmatrix} \ \ \ {\rm and} \ \ \ B=\begin{pmatrix}
1 & 1 & 1 & 1 & 1\\
1 & 1 & 1 & 1 & 1\\
0 & 0 & 0 & 0 & 0\\
1 & 1 & 1 & 1 & 1\\
0 & 0 & 0 & 0 & 0\\
\end{pmatrix}.
\]
Thus, the perturbed matrix $C$ and its Jordan canonical form $J$ (regardless of the permutations of its diagonal Jordan blocks) are
\[C=\begin{pmatrix}
2 & 2 & 1 & 1 & 1\\
1 & 2 & 2 & 1 & 1\\
0 & 0 & 1 & 0 & 0\\
1 & 1 & 1 & 2 & 1\\
0 & 0 & 0 & 0 & 1\\
\end{pmatrix} \ \ \ {\rm and} \ \ \ J=\begin{pmatrix}
1 & 1 & 0 & 0 & 0\\
0 & 1 & 0 & 0 & 0\\
0 & 0 & 1 & 0 & 0\\
0 & 0 & 0 & \frac{5-\sqrt{13}}{2} & 0\\
0 & 0 & 0 & 0 & \frac{5+\sqrt{13}}{2}\\
\end{pmatrix},\]
respectively. Since, $|\Lambda(A)|=1$, $d(A)=2$, ${\rm rank}(B)=1$, and $d(C)=1$, then the upper bounds in (1.1) and (1.2) are 
\[({\rm rank}(B)+1)|\Lambda(A)|+d(A)=4 \ \ \
{\rm and} 
\ \ \ ({\rm rank}(B)+1)|\Lambda(A)|+d(A)-d(C)=3,\] respectively. Notice that $|\Lambda(C)|=3$, which attains the upper bound in (1.2). 
 
\vskip 0.8cm

\noindent{\bf \large 2.\ Preliminaries}

\vskip 0.2cm

In this section, we introduce some basic notations and concepts which are frequently used in the subsequent content. Let $\lambda\in\Lambda(M)$. The multiplicity of $\lambda$ as a zero of the characteristic polynomial $f_{M}(\cdot)$ is called its \emph{algebraic multiplicity}. The dimension of the eigenspace of $M$ corresponding to $\lambda$ is called its \emph{geometric multiplicity}. Let $m_{a}(M,\lambda)$ and $m_{g}(M,\lambda)$ denote the algebraic and geometric multiplicity of $\lambda$ as an eigenvalue of $M$, respectively.

Actually, the geometric multiplicity of an eigenvalue $\lambda\in\Lambda(M)$ coincides with the number of Jordan blocks corresponding to $\lambda$ in the Jordan canonical form of $M$. Evidently, $m_{g}(M,\lambda)\geq1$ for all $\lambda\in\Lambda(M)$. If $m_{g}(M,\lambda)\equiv1$ for all $\lambda\in\Lambda(M)$, then $M$ is said to be \emph{nonderogatory}; otherwise, $M$ is called a \emph{derogatory} matrix. Recall that the geometric multiplicity of an eigenvalue is not greater than its algebraic multiplicity. If $m_{g}(M,\lambda)<m_{a}(M,\lambda)$ for some $\lambda\in\Lambda(M)$, then $M$ is called a \emph{defective} matrix. If $m_{g}(M,\lambda)=m_{a}(M,\lambda)$ for all $\lambda\in\Lambda(M)$, then $M$ is said to be \emph{nondefective}. Thus, a matrix $M$ is diagonalizable if and only if $M$ is nondefective; see, for example, [6, p. 58]. 

We next introduce the definitions of the \emph{defectivity} of an eigenvalue and a matrix; see [1, Definitions 1.1 and 1.2]. In addition, we define the \emph{derogatory index} of a matrix in Definition 2.3 below.

\vskip 0.2cm

\noindent{\bf Definition 2.1.} The \emph{defectivity of an eigenvalue} $\lambda\in\Lambda(M)$ is denoted by $d(M,\lambda)$, which is the difference between its algebraic and geometric multiplicities, i.e.,
\begin{align*}
d(M,\lambda):=m_{a}(M,\lambda)-m_{g}(M,\lambda).\tag{2.1}
\end{align*}

\vskip 0.2cm

\noindent{\bf Definition 2.2.} The \emph{defectivity of a matrix} $M$ is denoted by $d(M)$, which is the sum of the defectivities of its eigenvalues, i.e.,
\begin{align*}
d(M):=\sum_{\lambda\in\Lambda(M)}\left(m_{a}(M,\lambda)-m_{g}(M,\lambda)\right)=\sum_{\lambda\in\Lambda(M)}d(M,\lambda).\tag{2.2}
\end{align*}
Because of $m_{a}(M,\lambda)\geq m_{g}(M,\lambda)\geq1$ for all $\lambda\in\Lambda(M)$, we obtain $d(M,\lambda)\geq 0$ and $d(M)\geq 0$. Indeed, the defectivity of a matrix can be considered as a quantitative measure of its nondiagonalizability since $M$ is diagonalizable if and only if $d(M)=0$. Moreover, it is clear that the defectivity of a matrix is equal to the number of off-diagonal ones in its Jordan canonical form; see [1, Remark 1].

\vskip 0.2cm

\noindent{\bf Definition 2.3.} The \emph{derogatory index} of a matrix $M$ is denoted by $\mathcal{I}(M)$, which is defined by
\begin{align*}
\mathcal{I}(M):=\sum_{\lambda\in\Lambda(M)}\left(m_{g}(M,\lambda)-1\right).\tag{2.3}
\end{align*}
Let $M\in\mathbb{C}^{n\times n}$. Notice that since $\sum_{\lambda\in\Lambda(M)}m_{a}(M,\lambda)=n$, from (2.2), we have  $d(M)=n-\sum_{\lambda\in\Lambda(M)}m_{g}(M,\lambda)$. Thus, we obtain $|\Lambda(M)|+d(M)\leq n$ and $\mathcal{I}(M)=\sum_{\lambda\in\Lambda(M)}(m_{g}(M,\lambda)-1)=n-d(M)-|\Lambda(M)|$. Therefore, the quantity $\mathcal{I}(M)$ satisfies $0\leq\mathcal{I}(M)\leq n-1$ for all $M\in\mathbb{C}^{n\times n}$. We can find that $\mathcal{I}(M)=n-1$ if and only if $M$ is a scalar matrix (i.e., $M=cI$ for some $c\in\mathbb{C}$), and $\mathcal{I}(M)=0$ if and only if $M$ is nonderogatory.

We finally mention the relation between the degree of minimal polynomial and the number of distinct eigenvalues of a matrix. It is well known that a matrix $M$ can be diagonalized if and only if every zero of its minimal polynomial $q_{M}(t)$ has multiplicity one; see, for instance, [6, Corollary 3.3.10]. If $M$ is diagonalizable, then $|\Lambda(M)|$ is equal to the degree of $q_{M}(t)$ since all distinct eigenvalues of $M$ are the roots of $q_{M}(t)=0$; otherwise, $|\Lambda(M)|$ is strictly less than the degree of $q_{M}(t)$.
	
\vskip 0.8cm

\noindent{\bf \large 3.\ Main result}

\vskip 0.2cm

 We now give the main result of this paper, which provides an improved upper bound for the number of distinct eigenvalues of a perturbed matrix.
 
\vskip 0.2cm

\noindent{\bf Theorem 3.1.} \emph{Assume that $A,B\in\mathbb{C}^{n\times n}$ and let $C=A+B$, then} \[|\Lambda(C)|\leq({\rm rank}(B)+1)|\Lambda(A)|+d(A)-d(C).\] 

\noindent{\bf Proof.} We first define $S_{1}:=\Lambda(C)\cap\Lambda(A)$ and $S_{2}:=\Lambda(C)\backslash\Lambda(A)=\Lambda(C)\cap\left(\Lambda(A)\right)^{c}$, where $(\cdot)^{c}$ denotes the complement of a set. Then we have
\begin{align*}
|\Lambda(C)|=|\Lambda(C)\cap\Lambda(A)|+|\Lambda(C)\backslash\Lambda(A)|=|S_{1}|+|S_{2}|.\tag{3.1}
\end{align*} 

We next take into account the case that $S_{1}\neq\varnothing$ and $S_{2}\neq\varnothing$. Notice that an upper bound for $\sum_{\lambda\in S_{2}}m_{a}(C,\lambda)$ is established in the proof of Farrell’s work [1, Theorem 1.3] and it is crucial for Farrell's proof. Alternatively, we intend to seek an appropriate upper bound for the quantity $\sum_{\lambda\in S_{2}}m_{g}(C,\lambda)$. For any $M\in\mathbb{C}^{n\times n}$, we have $\sum_{\lambda\in\Lambda(M)}m_{a}(M,\lambda)=n$, and it follows that
\begin{align*}
\sum_{\lambda\in S_{1}}m_{a}(C,\lambda)+\sum_{\lambda\in S_{2}}m_{a}(C,\lambda)=\sum_{\lambda\in \Lambda(A)}m_{a}(C,\lambda)+\sum_{\lambda\in S_{2}}m_{a}(C,\lambda)=n,\tag{3.2}
\end{align*}
with the convention that $m_{a}(M,\lambda)=0$ (hence, $m_{g}(M,\lambda)=0$) if and only if $ \lambda\notin\Lambda(M)$. 

Let $\lambda$ be an arbitrary eigenvalue of $A$. If $\lambda\in S_{1}$, by ${\rm rank}(M_{1}+M_{2})\leq{\rm rank}(M_{1})+{\rm rank}(M_{2})$ for all $M_{1},M_{2}\in\mathbb{C}^{n\times n}$, we obtain
\[m_{g}(C,\lambda)=n-{\rm rank}(\lambda I-C)\geq n-{\rm rank}(\lambda I-A)-{\rm rank}(B)=m_{g}(A,\lambda)-{\rm rank}(B).\]
If $\lambda\in \Lambda(A)\backslash S_{1}$, then $\lambda I-C$ is nonsingular, which implies that $n-{\rm rank}(\lambda I-C)=0$, and  we have $m_{g}(C,\lambda)=0$ and $m_{g}(A,\lambda)-{\rm rank}(B)\leq0$.
Thus, we obtain
\begin{align*}
\sum_{\lambda\in \Lambda(A)}m_{a}(C,\lambda)&=\sum_{\lambda\in S_{1}}m_{a}(C,\lambda)\tag{3.3a}\\
&=\sum_{\lambda\in S_{1}}\left(m_{g}(C,\lambda)+d(C,\lambda)\right)\tag{3.3b}\\
&\geq\sum_{\lambda\in S_{1}}\left(m_{g}(A,\lambda)-{\rm 
rank}(B)+d(C,\lambda)\right)\\
&\quad+\sum_{\lambda\in \Lambda(A)\backslash S_{1}}\left(m_{g}(A,\lambda)-{\rm 
rank}(B)\right)\tag{3.3c}\\
&=\sum_{\lambda\in \Lambda(A)}\left(m_{g}(A,\lambda)-{\rm rank}(B)\right)+\sum_{\lambda\in S_{1}}d(C,\lambda)\tag{3.3d}\\
&=\sum_{\lambda\in \Lambda(A)}\left(m_{a}(A,\lambda)-{\rm rank}(B)-d(A,\lambda)\right)+\sum_{\lambda\in S_{1}}d(C,\lambda)\tag{3.3e}\\
&=n-|\Lambda(A)|\cdot{\rm rank}(B)-d(A)+\sum_{\lambda\in S_{1}}d(C,\lambda),\tag{3.3f}
\end{align*}
where we have used the definitions (2.1) and (2.2). It is worth mentioning that the key difference with the proof of Farrell [1] is the steps (3.3b) and (3.3c). By relation (3.2), we can obtain
\begin{align*}
\sum_{\lambda\in S_{2}}m_{g}(C,\lambda)+\sum_{\lambda\in S_{2}}d(C,\lambda)&=\sum_{\lambda\in S_{2}}m_{a}(C,\lambda)\tag{3.4a}\\
&\leq n-\left(n-|\Lambda(A)|\cdot{\rm rank}(B)-d(A)+\sum_{\lambda\in S_{1}}d(C,\lambda)\right)\tag{3.4b}\\
&=|\Lambda(A)|\cdot{\rm rank}(B)+d(A)-\sum_{\lambda\in S_{1}}d(C,\lambda).\tag{3.4c}
\end{align*}
Then
\[\sum_{\lambda\in S_{2}}m_{g}(C,\lambda)\leq|\Lambda(A)|\cdot{\rm rank}(B)+d(A)-\left(\sum_{\lambda\in S_{1}}d(C,\lambda)+\sum_{\lambda\in S_{2}}d(C,\lambda)\right),\]
in other words,
\begin{align*}
\sum_{\lambda\in S_{2}}m_{g}(C,\lambda)\leq|\Lambda(A)|\cdot{\rm rank}(B)+d(A)-d(C).\tag{3.5}
\end{align*}
Since $|S_{1}|\leq|\Lambda(A)|$ and $|S_{2}|\leq\sum_{\lambda\in S_{2}}m_{g}(C,\lambda)$, from (3.1) and (3.5), respectively, we deduce that
\begin{align*}
|\Lambda(C)|=|S_{1}|+|S_{2}|\leq({\rm rank}(B)+1)|\Lambda(A)|+d(A)-d(C).\tag{3.6}
\end{align*}

We finally consider two special cases: either $S_{1}$ or $S_{2}$ is empty. If $S_{1}=\varnothing$, i.e., $\Lambda(C)\cap\Lambda(A)=\varnothing$,  then $|S_{1}|=0$ and $S_{2}=\Lambda(C)$. Repeating the derivation in (3.4) yields
\[|\Lambda(A)|\cdot{\rm rank}(B)+d(A)\geq n.\]
Since
\[d(C)=\sum_{\lambda\in\Lambda(C)}d(C,\lambda)=\sum_{\lambda\in\Lambda(C)}\left(m_{a}(C,\lambda)-m_{g}(C,\lambda)\right)=n-\sum_{\lambda\in\Lambda(C)}m_{g}(C,\lambda)\]
and $m_{g}(C,\lambda)\geq1$ for all $\lambda\in\Lambda(C)$, we conclude that $|\Lambda(C)|+d(C)\leq n$ is always valid. Hence,
\begin{align*}
|\Lambda(C)|+d(C)\leq n<({\rm rank}(B)+1)|\Lambda(A)|+d(A).\tag{3.7}
\end{align*}
If $S_{2}=\varnothing$, i.e., $\Lambda(C)\subseteq\Lambda(A)$, then $|\Lambda(C)|\leq|\Lambda(A)|$, $S_{1}=\Lambda(C)$, and $|S_{2}|=0$. We repeat the derivation in (3.3) and obtain
\[\begin{split}
n&=\sum_{\lambda\in \Lambda(A)}m_{a}(C,\lambda)\\
&\geq\sum_{\lambda\in \Lambda(A)}\left(m_{a}(A,\lambda)-{\rm rank}(B)-d(A,\lambda)\right)+\sum_{\lambda\in \Lambda(C)}d(C,\lambda)\\
&=n-|\Lambda(A)|\cdot{\rm rank}(B) -d(A)+d(C),
\end{split}\]
then
\[d(C)\leq|\Lambda(A)|\cdot{\rm rank}(B)+d(A).\]
Using $|\Lambda(C)|\leq|\Lambda(A)|$ yields
\begin{align*}
|\Lambda(C)|+d(C)\leq({\rm rank}(B)+1)|\Lambda(A)|+d(A).\tag{3.8}
\end{align*}

Consequently, from (3.6)-(3.8), we conclude that
$|\Lambda(C)|\leq({\rm rank}(B)+1)|\Lambda(A)|+d(A)-d(C)$ always hold. This completes the proof. \ \ \ \ \ $\Box$

\vskip 0.2cm

\noindent{\bf Remark 3.2.} If the perturbed matrix $C$ is diagonalizable (i.e., $d(C)=0$), by Theorem 3.1, we obtain
\[|\Lambda(C)|\leq({\rm rank}(B)+1)|\Lambda(A)|+d(A),\]
which is exactly the inequality (1.1); otherwise, we have $d(C)\geq1$, then Theorem 3.1 gives
\[|\Lambda(C)|\leq({\rm rank}(B)+1)|\Lambda(A)|+d(A)-1,
\]
which is a smaller upper bound than the estimate (1.1).

\vskip 0.8cm

\noindent{\bf \large 4.\ Applications}

\vskip 0.2cm

In view of the improved upper bound for the number of distinct eigenvalues of a matrix after perturbation, we can establish some interesting results.

\vskip 0.2cm

The following Corollary 4.1 provides a lower bound for the derogatory index of a matrix after perturbation.

\vskip 0.2cm

\noindent{\bf Corollary 4.1.} \emph{Let $A,B\in\mathbb{C}^{n\times n}$ and $C=A+B$, the derogatory index of the perturbed matrix $C$ satisfies}
\[\mathcal{I}(C)\geq \mathcal{I}(A)-{\rm rank}(B)\cdot|\Lambda(A)|.\]

\noindent{\bf Proof.} We note that the definition (2.3) implies $\mathcal{I}(M)=n-(|\Lambda(M)|+d(M))$ for all $M\in\mathbb{C}^{n\times n}$. The statement follows immediately from Theorem 3.1. \ \ \ \ \ $\Box$

\vskip 0.2 cm

The next corollary plays an important role in the estimate for the number of Krylov iterations after a rank one update.

\vskip 0.2cm

\noindent{\bf Corollary 4.2.} \emph{Suppose that $A\in\mathbb{C}^{n\times n}$ is diagonalizable, ${\rm rank}(B)=1$, and let $C=A+B$. If $C$ is also diagonalizable, then $|\Lambda(C)|\leq2|\Lambda(A)|$. If $C$ is not diagonalizable, then  $|\Lambda(C)|\leq2|\Lambda(A)|-1$.}

\vskip 0.2cm

\noindent{\bf Proof.} On the basis of Theorem 3.1, we have $|\Lambda(C)|\leq2|\Lambda(A)|+d(A)-d(C)=2|\Lambda(A)|-d(C)$ since $A$ is diagonalizable. If $C$ is diagonalizable, then $d(C)=0$ which yields $|\Lambda(C)|\leq2|\Lambda(A)|$. If $C$ cannot be diagonalized, then $d(C)\geq1$ and we obtain $|\Lambda(C)|\leq2|\Lambda(A)|-1$.\ \ \ \ \ $\Box$

\vskip 0.2cm

\noindent{\bf Corollary 4.3.} \emph{Assume that $A,B\in\mathbb{C}^{n\times n}$ and let $C=A+B$. If $C$ is nonderogatory, then the number of distinct eigenvalues of $A$ satisfies} 
\[\frac{n-d(A)}{{\rm rank}(B)+1}\leq|\Lambda(A)|\leq n-d(A).\]

\noindent{\bf Proof.} Since $C$ is nonderogatory, we have $\mathcal{I}(C)=0$, i.e., $|\Lambda(C)|+d(C)=n$. According to Theorem 3.1, we obtain 
$n\leq({\rm rank}(B)+1)|\Lambda(A)|+d(A)$,
which implies that $|\Lambda(A)|\geq\frac{n-d(A)}{{\rm rank}(B)+1}$. Notice that $|\Lambda(A)|\leq n-d(A)$ is clear.\ \ \ \ \ $\Box$

\vskip 0.2cm

Given any $A\in\mathbb{C}^{n\times n}$, it can be decomposed as $A=\mathcal{H}(A)+\mathcal{S}(A)$, where $\mathcal{H}(A)=\frac{1}{2}(A+A^{\ast})$ and $\mathcal{S}(A)=\frac{1}{2}(A-A^{\ast})$ are the \emph{Hermitian} and \emph{skew-Hermitian} part of $A$, respectively. By Theorem 3.1, we can obtain the following inequality:

\vskip 0.2cm

\noindent{\bf Corollary 4.4.} \emph{Let $A\in\mathbb{C}^{n\times n}$, we have the following estimate}
\begin{align*}
|\Lambda(A)|\leq {\rm min}\big\{\left({\rm rank}\left(\mathcal{H}(A)\right)+1\right)\left|\Lambda\left(\mathcal{S}(A)\right)\right|, \left({\rm rank}\left(\mathcal{S}(A)\right)+1\right)\left|\Lambda\left(\mathcal{H}(A)\right)\right|\big\}-d(A).\tag{4.1}
\end{align*}

\noindent{\bf Proof.} Notice that both $\mathcal{H}(A)$ and $\mathcal{S}(A)$ are normal matrices, which can be unitarily diagonalized. Thus, the inequality (4.1) follows from Theorem 3.1. \ \ \ \ $\Box$

\vskip 0.2cm

\noindent{\bf Remark 4.5.} 
Since $\mathcal{H}(A)$ and $\mathcal{S}(A)$ are diagonalizable, we have $|\Lambda(\mathcal{H}(A))|\leq{\rm rank}(\mathcal{H}(A))+1$ and $|\Lambda(\mathcal{S}(A))|\leq{\rm rank}(\mathcal{S}(A))+1$. From inequality (4.1), we can deduce that
$|\Lambda(A)|\leq \left({\rm rank}\left(\mathcal{H}(A)\right)+1\right)\left({\rm rank}\left(\mathcal{S}(A)\right)+1\right)-d(A)$.

\vskip 0.2cm

\noindent{\bf Remark 4.6.} For any $A\in\mathbb{C}^{n\times n}$ and $\alpha\in\mathbb{C}$, $A$ can also be decomposed as \[A=(\mathcal{H}(A)+\alpha I)+(\mathcal{S}(A)-\alpha I)=\mathcal{H}_{\alpha}(A)+\mathcal{S}_{\alpha}(A),\]
where $\mathcal{H}_{\alpha}(A):=\mathcal{H}(A)+\alpha I$ and $\mathcal{S}_{\alpha}(A):=\mathcal{S}(A)-\alpha I$. Notice that $\mathcal{H}_{\alpha}(A)$ and $\mathcal{S}_{\alpha}(A)$ are normal for all $\alpha\in\mathbb{C}$, then $\mathcal{H}_{\alpha}(A)$ and $\mathcal{S}_{\alpha}(A)$ can be unitarily diagonalized. Since the parameter $\alpha$ is arbitrary, the estimate  (4.1) can be modified as 
\begin{align*}
|\Lambda(A)|\leq\inf_{\alpha\in\mathbb{C}}{\rm min}\big\{\left({\rm rank}\left(\mathcal{H}_{\alpha}(A)\right)+1\right)\left|\Lambda\left(\mathcal{S}_{\alpha}(A)\right)\right|, \left({\rm rank}\left(\mathcal{S}_{\alpha}(A)\right)+1\right)\left|\Lambda\left(\mathcal{H}_{\alpha}(A)\right)\right|\big\}-d(A).
\end{align*}
Similarly, we can obtain
\begin{align*}
|\Lambda(A)|\leq \inf_{\alpha\in\mathbb{C}}\big\{\left({\rm rank}\left(\mathcal{H}_{\alpha}(A)\right)+1\right)\left({\rm rank}\left(\mathcal{S}_{\alpha}(A)\right)+1\right)\big\}-d(A).
\end{align*}
Moreover, it is easy to see that other decompositions of $A$ can lead to some different estimates of $|\Lambda(A)|$.

\vskip 1.0cm

\noindent{\bf \large Acknowledgements}

\vskip 0.2cm

The author is indebted to Professor Farrell for his valuable comments and suggestions, which have greatly improved the original manuscript of this paper. The author would like to thank Professor Chen-Song Zhang for his helpful suggestions. This work was supported partially by the National Natural Science Foundation of China Grant 91430215 and 91530323.

\vskip 1.2cm

\noindent{\bf \large References}

\vskip 0.2cm

\small
{
\noindent{[1]}\ P. E. Farrell, The number of distinct eigenvalues of a matrix after perturbation, SIAM J. Matrix 

Anal. Appl., 37 (2016) 572–576.

\noindent{[2]}\ J. R. Bunch, C. P. Nielsen, D. C. Sorensen, Rank-one modification of the symmetric eigenproblem, 

Numer. Math., 31 (1978) 31–48.

\noindent{[3]}\ G. H. Golub, Some modified matrix eigenvalue problems, SIAM Rev., 15 (1973) 318–334.

\noindent{[4]}\ I. C. F. Ipsen, B. Nadler, Refined perturbation bounds for eigenvalues of Hermitian and non-Hermitian 

matrices, SIAM J. Matrix Anal. Appl., 31 (2009) 40–53.

\noindent{[5]}\ J. H. Wilkinson, The Algebraic Eigenvalue Problem, Monogr. Numer. Anal. 87, Oxford University 

Press, Oxford, 1965.

\noindent{[6]}\ R. A. Horn, C. R. Johnson, Matrix Analysis, Cambridge University Press, Cambridge, 1985.

}

\end{document}